# New Results for the Degree/Diameter Problem

Michael J. Dinneen
Computer Research and Applications
Los Alamos National Laboratory
Los Alamos, N.M. 87545
mjd@lanl.gov

Paul R. Hafner
Dept. of Mathematics and Statistics
University of Auckland
Auckland, New Zealand
hafner@mat.aukuni.ac.nz

**Abstract**

The results of computer searches for large graphs with given (small) degree and diameter are presented. The new graphs are Cayley graphs of semidirect products of cyclic groups and related groups. One fundamental use of our "dense graphs" is in the design of efficient communication network topologies.

## 1 Introduction

The design of interconnection networks comes up against two fundamental constraints: the number of connections which can be attached to any one node is limited (e.g. by lack of space) and the number of intermediate nodes on the communication route between two nodes must be small (e.g. to meet timing tolerances). While observing such constraints, one wishes to maximize the number of nodes which can participate in such a network. In the language of graph theory, this is the

**Degree/Diameter Problem:** *find graphs with maximal number of vertices with given constraints of maximum degree $\Delta$ and diameter $D$.*

All our graphs are undirected and we follow the definitions and notations of [11]. For convenience some basic definitions related to this paper are now presented.

The *degree* of a vertex is the number of incident edges to that vertex. The *distance* between two vertices is the length of a shortest path between them. The *diameter* of a graph is the maximum distance over all pairs of vertices. A graph with maximum degree $\Delta$ and diameter $D$ is called a $(\Delta, D)$-*graph*. The number of vertices in a graph is called the *order* of the graph.

The only known bound on the order $n$ of a graph with degree $\Delta$ and diameter $D$ is the *Moore bound* [4]:

$$n \leq 1 + \Delta + \Delta(\Delta - 1) + \cdots \Delta(\Delta - 1)^{D-1}$$
$$= \begin{cases} \dfrac{\Delta(\Delta - 1)^D - 2}{\Delta - 2} & \text{if } \Delta \neq 2, \\ 2D + 1 & \text{if } \Delta = 2. \end{cases}$$

Graphs whose order attains this bound are called *Moore graphs*.

It was proved in [24, 18, 5] that there do not exist any Moore graphs apart from the following:

| $\Delta$ | $D$ | $n$ | Description |
|---|---|---|---|
| 2 | $D$ | $2D+1$ | $(2D+1)$-gon |
| 3 | 2 | 10 | Petersen |
| 7 | 2 | 50 | Hoffman-Singleton |
| 57 | 2 | 3250 | ? |

It is not known if there exists a graph with $\Delta = 57$, $D = 2$ and $n = 3250$.

In the absence of realistic bounds, one resorts to constructions of largest possible $(\Delta, D)$-graphs. Descriptions of new constructions and corresponding listings of largest known $(\Delta, D)$-graphs have been published periodically by members of the mathematical and engineering communities, e.g. [26, 7, 20, 15, 8, 16, 9, 13, 12, 17, 23], beginning with a paper of Elspas [21]. Apart from the Moore graphs, the only graphs known to be optimal are $(3, 3)$, $(4, 2)$, and $(5, 2)$-graphs with 20, 15, and 24 vertices, respectively. This result goes back to Elspas [21], confirmed later by a general theorem [22, 6] which states that except for the case of a cycle on four vertices the number of vertices in a $(\Delta, D)$-graph never misses the Moore bound by one. In the current table of largest known graphs there is still a significant gap between the Moore bound and the orders of known graphs.

Some of the previously known largest $(\Delta, D)$-graphs were found as Cayley graphs, e.g. for classical groups [12, 16]. Others, e.g. cube-connected cycles, were quickly recognized as Cayley graphs [2]. According to a theorem of Sabidussi [25] every vertex-transitive graph can be realized as Cayley coset graph of some group. Cayley graphs have many desirable properties for engineering applications, e.g. vertex transitivity with associated ease of routing, and fault tolerance [1, 2, 15]

Considerable progress resulted from the first explicit use of semidirect products of cyclic groups [19]. These groups should be seen as capturing the essential features of the Borel groups used by [16] to establish their longstanding records for certain degree 4 graphs which are only now superseded by our results. The present paper continues this work and considers some additional types of groups. One of the major advantages of the semidirect products of cyclic groups is their abundance, together with the fact that they are easy to code on a computer.

The paper is organized as follows. In the next section we give some preliminary comments regarding Cayley graphs. The following section describes the types of groups that were used in finding dense Cayley graphs. We end the paper with a detailed listing of our new $(\Delta, D)$-graph records. For an overview of our results, see the "**DH**" entries in Table 1. This listing is maintained by J.-C. Bermond, C. Delorme and J.-J. Quisquater and available on request (e-mail: `cd@lri.fr`).

## 2 Cayley graphs

Let $G$ be a finite group, $S$ a subset of $G$ which generates $G$ and does not contain the identity. The *Cayley graph of $G$ with respect to $S$* is the directed graph whose vertices are the elements of $G$ and whose edge set is $E = \{(x,y): y = xs \text{ for some } s \in S\}$. If $S$ is closed under inverses, i.e. $S = S \cup S^{-1}$, then $(x,y) \in E$ if and only if $(y,x) \in E$. In this case the edges can be considered as undirected. In the present paper we deal only with this situation.

| $\Delta$ \ $D$ | 2 | 3 | 4 | 5 | 6 | 7 | 8 | 9 | 10 |
|---|---|---|---|---|---|---|---|---|---|
| 3 | P<br>10 | $C_5 * F_4$<br>20 | vC<br>38 | vC<br>70 | GFS<br>130 | $CR^*$<br>184 | $CR^*$<br>320 | 2cy<br>540 | 2cy<br>938 |
| 4 | $K_3 * C_5$<br>15 | Allwr<br>41 | $C_5 * C_{19}$<br>95 | $H'_3$<br>364 | $H_3(K_3)$<br>740 | **DH**<br>**1 155** | **DH**$^{**}$<br>**3 025** | **DH**<br>**7 550** | **DH**<br>**16 555** |
| 5 | $K_3 * X_8$<br>24 | Lente<br>70 | $Q_4(K_3)$<br>186 | $H'_3 d$<br>532 | $H_4(K_3)$<br>2 754 | **DH**<br>**5 334** | **DH**<br>**15 532** | **DH**<br>**49 932** | **DH**<br>**145 584** |
| 6 | $K_4 * X_8$<br>32 | $C_5 * C_{21}$<br>105 | **DH**$^*$<br>**360** | **DH**<br>**1 230** | $H_5(K_4)$<br>7 860 | **DH**<br>**18 775** | **DH**<br>**69 540** | **DH**<br>**275 540** | **DH**<br>**945 574** |
| 7 | HS<br>50 | **DH**$^*$<br>**144** | **DH**$^*$<br>**600** | **DH**<br>**2 756** | $H_4(K_4) < H_5$<br>10 566 | **DH**<br>**47 304** | **DH**<br>**214 500** | **DH**<br>**945 574** | Cam<br>4 773 696 |
| 8 | $P'_7$<br>57 | **DH**<br>**234** | **DH**<br>**1 012** | **DH**$^*$<br>**4 704** | $H_7(K_6)$<br>39 396 | **DH**<br>**127 134** | **DH**<br>**654 696** | **DH**$^{**}$<br>**2 408 704** | Cam<br>7 738 848 |
| 9 | $P'_8 d$<br>74 | $Q'_8$<br>585 | **DH**<br>**1 430** | **DH**<br>**7 344** | $H_8(K_6)$<br>75 198 | **DH**<br>**264 024** | **DH**$^{**}$<br>**1 354 896** | **DH**<br>**4 980 696** | Cam<br>19 845 936 |
| 10 | $P'_9$<br>91 | $Q'_8 d$<br>650 | **DH**<br>**2 200** | **DH**$^*$<br>**12 288** | $H_9(K_6)$<br>133 500 | **DH**<br>**554 580** | **DH**$^{**}$<br>**3 069 504** | **DH**<br>**9 003 000** | $Q_7 \Sigma_2 H_7$<br>47 059 200 |
| 11 | $P'_9 d$<br>94 | $Q'_8 d$<br>715 | $Q_7(T_4)$<br>3 200 | **DH**<br>**17 458** | $H_7(T_4)$<br>156 864 | **DH**<br>**945 574** | Cam<br>4 773 696 | Cam<br>25 048 800 | $Q_7 \Sigma_6 H_8$<br>179 755 200 |
| 12 | $P'_{11}$<br>133 | $Q'_8 d$<br>780 | $Q'_8 * X_8$<br>4 680 | **DH**<br>**26 871** | $H_{11}(K_6)$<br>355 812 | Dinn<br>1 732 514 | **DH**<br>**10 007 820** | **DH**<br>**48 532 122** | $Q_8 \Sigma_6 H_9$<br>466 338 600 |
| 13 | $P'_{11} d$<br>136 | $Q'_8 d$<br>845 | $Q_9(T_4)$<br>6 560 | **DH**<br>**37 056** | $H_9(T_4)$<br>531 440 | Cam<br>2 723 040 | **DH**<br>**15 027 252** | **DH**<br>**72 598 920** | $Q_9 \Sigma_6 H_9$<br>762 616 400 |
| 14 | $P'_{13}$<br>183 | $Q'_8 d$<br>910 | $Q_9(T_5)$<br>8 200 | **DH**<br>**53 955** | $H_{13}(K_7)$<br>806 636 | $K_1 \Sigma_8 H_{11}$<br>6 200 460 | Dinn<br>29 992 052 | $P_9 \Sigma_7 H_{11}$<br>164 755 080 | $Q_8 \Sigma_6 H_{11}$<br>1 865 452 680 |
| 15 | $P'_{13} d$<br>186 | $(\otimes Q_{2,4})'$<br>1 215 | $Q_{11}(T_4)$<br>11 712 | **DH**<br>**69 972** | $H_{11}(T_4)$<br>1 417 248 | **DH**<br>**7 100 796** | **DH**<br>**38 471 006** | $P_{11} \Sigma_1 H_{11}$<br>282 740 976 | $Q_{11} \Sigma_6 H_{11}$<br>3 630 989 376 |
| 16 | $P'_{13} d$<br>197 | $(\otimes Q_3)'$<br>1 600 | $Q_{11}(T_5)$<br>14 640 | $(\otimes H_3)'$<br>132 496 | $H_{11}(T_5)$<br>1 771 560 | $K_1 \Sigma_8 H_{13}$<br>14 882 658 | $K_{9,9} \Sigma_6 H_{13}$<br>86 882 544 | $P_9 \Sigma_7 H_{11}$<br>585 652 704 | $Q_{11} d \Sigma_6 H_{13}$<br>7 394 669 856 |

Table 1: Table of largest known ($\Delta$, $D$)-graphs

An undirected Cayley graph is $\Delta$-regular (i.e. each vertex has degree $\Delta$) where $\Delta$ is the number of generators in $S = S \cup S^{-1}$. Since Cayley graphs are obviously vertex transitive [25], the computation of the diameter $D$ requires only the determination of the distances from the identity element to all other elements.

It is obvious that for the generation of large $(\Delta, D)$-graphs one has to avoid abelian groups: the relation $ab = ba$ prevents the desired Moore bound even at distance two. It is rather surprising that good results can be achieved with the following nonabelian groups which are easily constructed from cyclic (abelian) groups.

# 3 Description of the groups

## 3.1 Semidirect product of $Z_m$ with $Z_n$

Let $Z_n$ be the cyclic group of integers under addition modulo $n$. Every automorphism of $Z_n$ can be represented by a unit of the ring $Z_n$. If the multiplicative order of the unit $a \in Z_n$ divides $m$ a semidirect product of $Z_m$ with $Z_n$ can be defined by

$$[x,y][u,v] = [x + u \bmod m, ya^u + v \bmod n].$$

In Table 1 groups of this kind are identified by **DH**; our detailed listing uses the symbol $m \times_a n$, e.g. $52 \times_2 53$.

These groups can be described in terms of generators and relations in the following way:

$$m \times_a n = \langle\, x, y \mid x^m = y^n = x^{-1}yxy^{-a} = 1 \,\rangle. \tag{1}$$

## 3.2 Semidirect product of $Z_m$ with $Z_n \times Z_n$

A construction which brought several improvements at the low end of the table are semidirect products of a cyclic group $Z_m$ with a direct sum $Z_n \times Z_n$. An automorphism $\sigma$ of $Z_n \times Z_n$ is determined by the images of the generators $\sigma([1,0]) = [x,y]$ and $\sigma([0,1]) = [z,t]$. If the order of $\sigma$ divides $m$ we can define a multiplication on $Z_m \times Z_n \times Z_n$ by:

$$[c,d,e][f,g,h] = [c + f \bmod m,\ [d,e]\begin{bmatrix} x & y \\ z & t \end{bmatrix}^f + [g,h] \bmod n].$$

In Table 1 groups of this kind are identified by **DH**$^*$; our detailed listing uses the symbol $m \times_\sigma n^2$, e.g. $40 \times_\sigma 3^2$; the action of the cyclic group is not encoded into this symbol and is specified separately.

These groups can be described in terms of generators and relations in the following way ($x$, $y$, $z$, $t$ have the same meaning as before):

$$m \times_\sigma n^2 = \langle\, a, b, c \mid a^m = b^n = c^n = bcb^{-1}c^{-1} = a^{-1}bac^{-y}b^{-x} = a^{-1}cac^{-t}b^{-z} = 1\,\rangle. \qquad (2)$$

## 3.3 Semidirect product of $Z_m \times_a Z_n$ with itself

The next step in these explorations was to double up on our first construction: let $G = m \times_a n$ and consider the semidirect product $G \times_\sigma G$ where $G$ acts on itself by conjugation. Explicitly, the multiplication of quadruples is given by

$$[x_1, x_2, x_3, x_4][y_1, y_2, y_3, y_4] = \big[[x_1, x_2][y_1, y_2],\, [y_1, y_2]^{-1}[x_3, x_4][y_1, y_2][y_3, y_4]\big] =$$
$$= [x_1 + y_1,\, x_2 a^{y_1} + y_2,\, x_3 + y_3,\, x_4 a^{y_1 + y_3} + y_2 a^{y_3} - y_2 a^{x_3 + y_3} + y_4].$$

Multiplication of pairs is the multiplication in $m \times_a n$. In Table 1 groups of this kind are identified by **DH**$^{**}$; our detailed listing uses the symbol $[m \times_a n]^2_\sigma$, e.g. $[5 \times_4 11]^2_\sigma$.

These groups can be described in terms of generators and relations in the following way:

$$\begin{aligned}[m \times_a n]^2_\sigma &= \langle\, r, s, t, u \mid r^m = s^n = r^{-1}srs^{-a} = t^m = u^n = t^{-1}utu^{-a} = \qquad (3)\\ &= r^{-1}trt^{-1} = s^{-1}tsu^{-1}t^{-1}u = r^{-1}uru^{-a} = s^{-1}usu^{-1} = 1\,\rangle.\end{aligned}$$

# 4 Final Remarks

Approximately two to four days of computer time was consumed for each $(\Delta, D)$ entry of Table 1. Primarily Sun-4s were used in our search for dense graphs with one million nodes or less while Los Alamos National Laboratory's Cray Y-MPs tackled most of the diameter computations for the larger Cayley graphs. For a fixed target group, several thousand sets of random generators were tried – with fewer combinations checked as the order of the groups grew.

Not all of our attempts were crowned with success. We note the following two types of groups which did not stand a chance in the competition.

The full automorphism group of the cyclic group $Z_n$ is isomorphic to the group $U_n$ of units modulo $n$. Like every finite abelian group, $U_n$ can be decomposed:

$$U_n = C_1 \oplus C_2 \oplus \cdots \oplus C_k,$$

where each $C_i$ is cyclic of order $d_i$ and $d_1|d_2|\cdots|d_k$. The groups $m \times_a n$ were concerned with $C_k$ or a subgroup thereof. Letting $C_{k-1} \oplus C_k$ act on $Z_n$ did not lead to any large graphs.

Another natural step after the groups $m \times_\sigma n^2$ would be a cyclic group $Z_m$ operating on $Z_n \times Z_n \times Z_n$. Again, these constructions did not produce good graphs.

# 5 Acknowledgments

The second author is indebted to Marston Conder for introducing him to the problem. In particular, the suggestion of looking at groups $m \times_\sigma n^2$ is due to him. Peter Dobcsányi provided help with some of the subtleties of C programming. The computer algebra package CAYLEY [14] was used to prepare some ancillary files and to check the results for some of the smaller groups.

# A  Groups and Generators

For an explanation of the entries in the column headed 'Group' refer to equations (1), (2), (3) in section 3.

| $(\Delta, D)$ | Order | Group | Generators | Inverses | Order of Generator |
|---|---|---|---|---|---|
| ( 4, 7 ) | 1155 | $15 \times_4 77$ | [ 6 2 ] | [ 9 5 ] | 35 |
| | | | [ 10 9 ] | [ 5 24 ] | 33 |
| ( 4, 8 ) | 3025 | $[5 \times_4 11]_\sigma^2$ | [ 0 5 1 1 ] | [ 0 6 4 7 ] | 55 |
| | | | [ 4 5 1 7 ] | [ 1 2 4 0 ] | 55 |
| ( 4, 9 ) | 7550 | $25 \times_{171} 302$ | [ 8 156 ] | [ 17 82 ] | 25 |
| | | | [ 10 31 ] | [ 15 285 ] | 10 |
| ( 4, 10 ) | 16555 | $35 \times_{256} 473$ | [ 8 342 ] | [ 27 343 ] | 35 |
| | | | [ 3 60 ] | [ 32 373 ] | 35 |
| ( 4, 11 ) | 42861 | $39 \times_{16} 1099$ | [ 19 863 ] | [ 20 783 ] | 39 |
| | | | [ 28 466 ] | [ 11 544 ] | 39 |
| ( 4, 12 ) | 95634 | $66 \times_2 1449$ | [ 16 1289 ] | [ 50 1270 ] | 99 |
| | | | [ 13 1253 ] | [ 53 854 ] | 66 |
| ( 4, 13 ) | 140868 | $84 \times_2 1677$ | [ 37 462 ] | [ 47 1212 ] | 84 |
| | | | [ 44 814 ] | [ 40 119 ] | 21 |
| ( 5, 7 ) | 5334 | $42 \times_{27} 127$ | [ 8 50 ] | [ 34 10 ] | 21 |
| | | | [ 27 15 ] | [ 15 99 ] | 14 |
| | | | [ 21 0 ] | | 2 |
| ( 5, 8 ) | 15532 | $44 \times_{207} 353$ | [ 25 50 ] | [ 19 200 ] | 44 |
| | | | [ 43 41 ] | [ 1 338 ] | 44 |
| | | | [ 22 0 ] | | 2 |
| ( 5, 9 ) | 49932 | $36 \times_{25} 1387$ | [ 7 440 ] | [ 29 462 ] | 36 |
| | | | [ 28 105 ] | [ 8 1113 ] | 9 |
| | | | [ 18 0 ] | | 2 |
| ( 5, 10 ) | 145584 | $48 \times_{2300} 3033$ | [ 6 2477 ] | [ 42 1303 ] | 72 |
| | | | [ 17 2951 ] | [ 31 2786 ] | 48 |
| | | | [ 24 0 ] | | 2 |
| ( 6, 4 ) | 360 | $40 \times_\sigma 3^2$ | [ 39 1 2 ] | [ 1 0 2 ] | 40 |
| | | [ 1 0 ] $\to$ [ 1 1 ] | [ 14 2 0 ] | [ 26 2 1 ] | 20 |
| | | [ 0 1 ] $\to$ [ 1 0 ] | [ 32 0 2 ] | [ 8 0 1 ] | 15 |
| ( 6, 5 ) | 1230 | $15 \times_{37} 82$ | [ 7 60 ] | [ 8 68 ] | 15 |
| | | | [ 4 30 ] | [ 11 38 ] | 15 |
| | | | [ 3 37 ] | [ 12 23 ] | 10 |
| ( 6, 7 ) | 18775 | $25 \times_{481} 751$ | [ 9 556 ] | [ 16 291 ] | 25 |
| | | | [ 2 570 ] | [ 23 386 ] | 25 |
| | | | [ 3 22 ] | [ 22 605 ] | 25 |

| $(\Delta, D)$ | Order | Group | Generators | Inverses | Order of Generator |
|---|---|---|---|---|---|
| ( 6, 8 ) | 69540 | $60 \times_8 1159$ | [ 42 500 ] | [ 18 545 ] | 190 |
| | | | [ 23 1038 ] | [ 37 94 ] | 60 |
| | | | [ 57 403 ] | [ 3 1125 ] | 20 |
| ( 6, 9 ) | 275540 | $92 \times_{2202} 2995$ | [ 28 2233 ] | [ 64 2932 ] | 115 |
| | | | [ 51 2790 ] | [ 41 325 ] | 92 |
| | | | [ 42 2831 ] | [ 50 556 ] | 46 |
| ( 6, 10 ) | 945574 | $238 \times_{81} 3973$ | [ 211 2137 ] | [ 27 216 ] | 238 |
| | | | [ 30 32 ] | [ 208 2390 ] | 119 |
| | | | [ 4 460 ] | [ 234 412 ] | 119 |
| ( 7, 3 ) | 144 | $16 \times_\sigma 3^2$ | [ 13 0 1 ] | [ 3 1 2 ] | 16 |
| | | $[ 1\ 0 ] \to [ 1\ 1 ]$ | [ 5 2 2 ] | [ 11 2 0 ] | 16 |
| | | $[ 0\ 1 ] \to [ 1\ 0 ]$ | [ 10 0 2 ] | [ 6 2 2 ] | 8 |
| | | | [ 8 0 0 ] | | 2 |
| ( 7, 4 ) | 600 | $24 \times_\sigma 5^2$ | [ 22 0 3 ] | [ 2 3 1 ] | 12 |
| | | $[ 1\ 0 ] \to [ 1\ 2 ]$ | [ 15 1 3 ] | [ 9 0 1 ] | 8 |
| | | $[ 0\ 1 ] \to [ 4\ 0 ]$ | [ 18 4 2 ] | [ 6 2 1 ] | 4 |
| | | | [ 12 0 0 ] | | 2 |
| ( 7, 5 ) | 2756 | $52 \times_2 53$ | [ 25 45 ] | [ 27 37 ] | 52 |
| | | | [ 30 23 ] | [ 22 18 ] | 26 |
| | | | [ 40 39 ] | [ 12 51 ] | 13 |
| | | | [ 26 0 ] | | 2 |
| ( 7, 7 ) | 47304 | $72 \times_5 657$ | [ 67 155 ] | [ 5 491 ] | 72 |
| | | | [ 59 160 ] | [ 13 100 ] | 72 |
| | | | [ 66 305 ] | [ 6 253 ] | 36 |
| | | | [ 36 0 ] | | 2 |
| ( 7, 8 ) | 214500 | $60 \times_2 3575$ | [ 50 1706 ] | [ 10 1231 ] | 66 |
| | | | [ 29 3164 ] | [ 31 2428 ] | 60 |
| | | | [ 56 3360 ] | [ 4 3440 ] | 15 |
| | | | [ 30 0 ] | | 2 |
| ( 7, 9 ) | 945574 | $238 \times_{81} 3973$ | [ 71 3406 ] | [ 167 1196 ] | 238 |
| | | | [ 109 2984 ] | [ 129 2397 ] | 238 |
| | | | [ 184 915 ] | [ 54 441 ] | 119 |
| | | | [ 119 0 ] | | 2 |
| ( 8, 3 ) | 234 | $18 \times_3 13$ | [ 7 5 ] | [ 11 7 ] | 18 |
| | | | [ 5 1 ] | [ 13 10 ] | 18 |
| | | | [ 16 12 ] | [ 2 9 ] | 9 |
| | | | [ 14 2 ] | [ 4 7 ] | 9 |
| ( 8, 4 ) | 1012 | $22 \times_{25} 46$ | [ 1 7 ] | [ 21 31 ] | 22 |
| | | | [ 14 33 ] | [ 8 39 ] | 22 |
| | | | [ 18 19 ] | [ 4 41 ] | 22 |
| | | | [ 4 44 ] | [ 18 26 ] | 11 |

| ($\Delta$, D) | Order | Group | Generators | Inverses | Order of Generator |
|---|---|---|---|---|---|
| ( 8, 5 ) | 4704 | $96 \times_\sigma 7^2$ | [ 5 3 2 ] | [ 91 1 5 ] | 96 |
| | | [ 1 0 ] $\to$ [ 1 6 ] | [ 1 5 3 ] | [ 95 6 2 ] | 96 |
| | | [ 0 1 ] $\to$ [ 5 5 ] | [ 39 0 4 ] | [ 57 3 3 ] | 32 |
| | | | [ 36 6 1 ] | [ 60 0 4 ] | 8 |
| ( 8, 7 ) | 127134 | $126 \times_{993} 1009$ | [ 73 719 ] | [ 53 953 ] | 126 |
| | | | [ 60 639 ] | [ 66 1007 ] | 21 |
| | | | [ 48 998 ] | [ 78 745 ] | 21 |
| | | | [ 14 447 ] | [ 112 711 ] | 9 |
| ( 8, 8 ) | 654696 | $216 \times_{625} 3031$ | [ 160 1966 ] | [ 56 816 ] | 27 |
| | | | [ 14 2452 ] | [ 202 885 ] | 108 |
| | | | [ 127 1541 ] | [ 89 1305 ] | 216 |
| | | | [ 153 1702 ] | [ 63 678 ] | 168 |
| ( 8, 9 ) | 2408704 | $[16 \times_8 97]^2_\sigma$ | [ 10 59 6 89 ] | [ 6 57 10 27 ] | 776 |
| | | | [ 3 7 14 92 ] | [ 13 68 2 21 ] | 16 |
| | | | [ 1 41 11 79 ] | [ 15 7 5 44 ] | 16 |
| | | | [ 5 80 6 13 ] | [ 11 26 10 81 ] | 16 |
| ( 9, 4 ) | 1430 | $10 \times_{64} 143$ | [ 0 84 ] | [ 0 59 ] | 143 |
| | | | [ 7 54 ] | [ 3 80 ] | 10 |
| | | | [ 1 51 ] | [ 9 64 ] | 10 |
| | | | [ 7 121 ] | [ 3 121 ] | 10 |
| | | | [ 5 0 ] | | 2 |
| ( 9, 5 ) | 7344 | $48 \times_5 153$ | [ 0 71 ] | [ 0 82 ] | 153 |
| | | | [ 16 86 ] | [ 32 118 ] | 153 |
| | | | [ 47 97 ] | [ 1 127 ] | 48 |
| | | | [ 37 130 ] | [ 11 100 ] | 48 |
| | | | [ 24 0 ] | | 2 |
| ( 9, 7 ) | 264024 | $72 \times_{1923} 3667$ | [ 27 3187 ] | [ 45 2969 ] | 152 |
| | | | [ 1 1495 ] | [ 71 1151 ] | 72 |
| | | | [ 7 1659 ] | [ 65 2792 ] | 72 |
| | | | [ 6 1431 ] | [ 66 661 ] | 12 |
| | | | [ 36 0 ] | | 2 |
| ( 9, 8 ) | 1354896 | $[12 \times_6 97]^2_\sigma$ | [ 11 34 11 21 ] | [ 1 87 1 70 ] | 12 |
| | | | [ 5 75 0 39 ] | [ 7 62 0 40 ] | 12 |
| | | | [ 8 76 3 60 ] | [ 4 56 9 55 ] | 12 |
| | | | [ 6 22 10 48 ] | [ 6 22 2 85 ] | 6 |
| | | | [ 6 0 0 0 ] | | 2 |
| ( 9, 9 ) | 4980696 | $1288 \times_{11} 3867$ | [ 442 2170 ] | [ 846 2609 ] | 1932 |
| | | | [ 925 2708 ] | [ 363 3857 ] | 1288 |
| | | | [ 1276 3002 ] | [ 12 2002 ] | 966 |
| | | | [ 408 2678 ] | [ 880 619 ] | 483 |
| | | | [ 644 0 ] | | 2 |

| $(\Delta, D)$ | Order | Group | Generators | Inverses | Order of Generator |
|---|---|---|---|---|---|
| ( 10, 4 ) | 2200 | $20 \times_3 110$ | [ 17 1 ] | [ 3 83 ] | 20 |
| | | | [ 3 0 ] | [ 17 0 ] | 20 |
| | | | [ 2 80 ] | [ 18 40 ] | 10 |
| | | | [ 4 33 ] | [ 16 77 ] | 10 |
| | | | [ 8 6 ] | [ 12 34 ] | 5 |
| ( 10, 5 ) | 12288 | $48 \times_\sigma 16^2$ | [ 25 7 15 ] | [ 23 15 6 ] | 48 |
| | | $[ 1\ 0 ] \to [ 1\ 15 ]$ | [ 29 1 2 ] | [ 19 13 5 ] | 48 |
| | | $[ 0\ 1 ] \to [ 7\ 8 ]$ | [ 43 0 9 ] | [ 5 5 5 ] | 48 |
| | | | [ 46 6 10 ] | [ 2 14 12 ] | 24 |
| | | | [ 44 10 15 ] | [ 4 11 8 ] | 12 |
| ( 10, 7 ) | 554580 | $156 \times_2 3555$ | [ 32 2159 ] | [ 124 1096 ] | 585 |
| | | | [ 66 2090 ] | [ 90 3355 ] | 234 |
| | | | [ 4 1182 ] | [ 152 2148 ] | 195 |
| | | | [ 41 2287 ] | [ 115 3319 ] | 156 |
| | | | [ 3 3039 ] | [ 153 1842 ] | 52 |
| ( 10, 8 ) | 3069504 | $[24 \times_{52} 73]^2_\sigma$ | [ 19 36 22 45 ] | [ 5 15 2 8 ] | 24 |
| | | | [ 3 47 7 41 ] | [ 21 61 17 43 ] | 24 |
| | | | [ 19 61 4 69 ] | [ 5 68 20 17 ] | 24 |
| | | | [ 14 37 22 39 ] | [ 10 61 2 63 ] | 12 |
| | | | [ 18 3 10 26 ] | [ 6 65 14 30 ] | 12 |
| ( 10, 9 ) | 9003000 | $3000 \times_{14} 3001$ | [ 77 967 ] | [ 2923 395 ] | 3000 |
| | | | [ 1864 494 ] | [ 1136 124 ] | 375 |
| | | | [ 1624 838 ] | [ 1376 1044 ] | 375 |
| | | | [ 2380 572 ] | [ 620 2799 ] | 150 |
| | | | [ 576 73 ] | [ 2424 1225 ] | 125 |
| ( 11, 5 ) | 17458 | $14 \times_{729} 1247$ | [ 1 459 ] | [ 13 1154 ] | 14 |
| | | | [ 1 134 ] | [ 13 1228 ] | 14 |
| | | | [ 3 433 ] | [ 11 1199 ] | 14 |
| | | | [ 10 443 ] | [ 4 910 ] | 7 |
| | | | [ 10 325 ] | [ 4 175 ] | 7 |
| | | | [ 7 0 ] | | 2 |
| ( 11, 7 ) | 945574 | $238 \times_{81} 3973$ | [ 111 2465 ] | [ 127 2668 ] | 238 |
| | | | [ 131 211 ] | [ 107 3540 ] | 238 |
| | | | [ 59 3508 ] | [ 179 513 ] | 238 |
| | | | [ 32 3841 ] | [ 206 1445 ] | 119 |
| | | | [ 188 2240 ] | [ 50 3522 ] | 119 |
| | | | [ 119 0 ] | | 2 |

| ($\Delta$, D) | Order | Group | Generators | Inverses | Order of Generator |
|---|---|---|---|---|---|
| ( 12, 5 ) | 26871 | $39 \times_{16} 689$ | [ 13 383 ] | [ 26 518 ] | 159 |
| | | | [ 5 667 ] | [ 34 235 ] | 39 |
| | | | [ 28 303 ] | [ 11 16 ] | 39 |
| | | | [ 25 41 ] | [ 14 86 ] | 39 |
| | | | [ 36 361 ] | [ 3 627 ] | 13 |
| | | | [ 27 400 ] | [ 12 81 ] | 13 |
| ( 12, 8 ) | 10007820 | $2580 \times_5 3879$ | [ 428 151 ] | [ 2152 2690 ] | 1935 |
| | | | [ 678 3117 ] | [ 1902 3102 ] | 1290 |
| | | | [ 276 717 ] | [ 2304 3621 ] | 645 |
| | | | [ 1375 2266 ] | [ 1205 2425 ] | 516 |
| | | | [ 735 3686 ] | [ 1845 3830 ] | 172 |
| | | | [ 1665 1590 ] | [ 915 2058 ] | 172 |
| ( 12, 9 ) | 48532122 | $6966 \times_5 6967$ | [ 4643 2316 ] | [ 2323 4128 ] | 6966 |
| | | | [ 5284 3171 ] | [ 1682 2937 ] | 3483 |
| | | | [ 69 877 ] | [ 6897 6154 ] | 2322 |
| | | | [ 1812 5557 ] | [ 5154 3426 ] | 1161 |
| | | | [ 3987 6878 ] | [ 2979 1122 ] | 774 |
| | | | [ 6880 2878 ] | [ 86 54 ] | 81 |
| ( 13, 5 ) | 37056 | $64 \times_{125} 579$ | [ 23 254 ] | [ 41 161 ] | 64 |
| | | | [ 27 422 ] | [ 37 44 ] | 64 |
| | | | [ 9 486 ] | [ 55 456 ] | 64 |
| | | | [ 43 278 ] | [ 21 116 ] | 64 |
| | | | [ 8 124 ] | [ 56 428 ] | 24 |
| | | | [ 12 156 ] | [ 52 159 ] | 16 |
| | | | [ 32 0 ] | | 2 |
| ( 13, 8 ) | 15027252 | $3876 \times_2 3877$ | [ 2785 2526 ] | [ 1091 3330 ] | 3876 |
| | | | [ 3739 2740 ] | [ 137 1997 ] | 3876 |
| | | | [ 790 430 ] | [ 3086 1901 ] | 1938 |
| | | | [ 2560 3681 ] | [ 1316 3504 ] | 969 |
| | | | [ 1716 2226 ] | [ 2160 3789 ] | 323 |
| | | | [ 1520 223 ] | [ 2356 1970 ] | 51 |
| | | | [ 1938 0 ] | | 2 |
| ( 13, 9 ) | 72598920 | $8520 \times_{13} 8521$ | [ 110 233 ] | [ 8410 7302 ] | 852 |
| | | | [ 294 989 ] | [ 8226 4335 ] | 1420 |
| | | | [ 1492 3266 ] | [ 7028 5171 ] | 2130 |
| | | | [ 4967 2870 ] | [ 3553 5684 ] | 8520 |
| | | | [ 51 1308 ] | [ 8469 361 ] | 2840 |
| | | | [ 6648 4785 ] | [ 1872 8005 ] | 355 |
| | | | [ 4260 0 ] | | 2 |

| (Δ, D) | Order | Group | Generators | Inverses | Order of Generator |
|---|---|---|---|---|---|
| ( 14, 5 ) | 53955 | $45 \times_{234} 1199$ | [ 36 972 ] | [ 9 336 ] | 545 |
| | | | [ 2 533 ] | [ 43 1119 ] | 45 |
| | | | [ 32 288 ] | [ 13 617 ] | 45 |
| | | | [ 4 529 ] | [ 41 481 ] | 45 |
| | | | [ 1 508 ] | [ 44 377 ] | 45 |
| | | | [ 34 850 ] | [ 11 618 ] | 45 |
| | | | [ 16 1187 ] | [ 29 565 ] | 45 |
| ( 15, 5 ) | 69972 | $84 \times_{81} 833$ | [ 80 62 ] | [ 4 312 ] | 357 |
| | | | [ 40 365 ] | [ 44 502 ] | 357 |
| | | | [ 3 565 ] | [ 81 492 ] | 196 |
| | | | [ 45 636 ] | [ 39 652 ] | 196 |
| | | | [ 15 453 ] | [ 69 520 ] | 196 |
| | | | [ 50 644 ] | [ 34 168 ] | 42 |
| | | | [ 35 219 ] | [ 49 468 ] | 12 |
| | | | [ 42 0 ] | | 2 |
| ( 15, 7 ) | 7100796 | $1884 \times_{49} 3769$ | [ 943 2746 ] | [ 941 1979 ] | 1884 |
| | | | [ 1402 2317 ] | [ 482 2419 ] | 942 |
| | | | [ 34 1747 ] | [ 1850 2321 ] | 942 |
| | | | [ 1131 40 ] | [ 753 3119 ] | 628 |
| | | | [ 1040 3137 ] | [ 844 2872 ] | 471 |
| | | | [ 1568 145 ] | [ 316 2312 ] | 471 |
| | | | [ 342 393 ] | [ 1542 461 ] | 314 |
| | | | [ 942 0 ] | | 2 |
| ( 15, 8 ) | 38471006 | $6202 \times_2 6203$ | [ 4159 3195 ] | [ 2043 3638 ] | 6202 |
| | | | [ 1486 1096 ] | [ 4716 3827 ] | 3101 |
| | | | [ 3916 3469 ] | [ 2286 3796 ] | 3101 |
| | | | [ 4333 2871 ] | [ 1869 276 ] | 886 |
| | | | [ 5215 3672 ] | [ 987 4456 ] | 886 |
| | | | [ 2968 321 ] | [ 3234 5123 ] | 443 |
| | | | [ 5726 2428 ] | [ 476 4180 ] | 443 |
| | | | [ 3101 0 ] | | 2 |

# References

[1] S. B. Akers and B. Krishnamurthy, "On group graphs and their fault-tolerance," *IEEE Trans. Computers*, **36** (1987), 885–888.

[2] S. B. Akers and B. Krishnamurthy, "A group-theoretic model for symmetric interconnection networks," *IEEE Trans. Computers*, **38** (1989), 555–565.

[3] R. Bar-Yehuda and T. Etzion, "Connections between two cycles – a new design of dense processor interconnection networks," *Discrete Applied Math.*, **37/38** (1992), 29–43.

[4] N. Biggs, *Algebraic Graph Theory*, Cambridge University Press, Cambridge (1974).

[5] E. Bannai and T. Ito, "On finite Moore graphs," *J. Fac. Sci. Univ. Tokyo*, **20** (1973), 191–208.

[6] E. Bannai and T. Ito, "Regular graphs with excess one," *Discrete Math.*, **37** (1981), 147–158.

[7] J.-C. Bermond, C. Delorme and J.-J. Quisquater, "Tables of large graphs with given degree and diameter," *Information Processing Letters*, **15** (1982), 10–13.

[8] J.-C. Bermond, C. Delorme and J.-J. Quisquater, "Strategies for interconnection networks: some methods from graph theory," *Journal of Parallel and Distributed Computing*, **3** (1986), 433–449.

[9] J.-C. Bermond, C. Delorme and J.-J. Quisquater, "Table of large $(\Delta, D)$-graphs," *Discrete Applied Math.*, **37/38** (1992), 575–577.

[10] J. Bond, C. Delorme and W.F. de La Vega, "Large Cayley graphs with small degree and diameter," *Rapport de Recherche no. 392*, LRI, Orsay, 1987.

[11] J. Bondy and U. Murty, *Graph Theory with Applications*, American Elsevier Publishing Co., 1976.

[12] L. Campbell, G. E. Carlsson, M. J. Dinneen, V. Faber, M. R. Fellows, M. A. Langston, J. W. Moore, A. P. Mullhaupt and H. B. Sexton, "Small diameter symmetric networks from linear groups," *IEEE Trans. Computers*, **41** (1992) 218-220.

[13] L. Campbell, "Dense group networks," *Discrete Applied Math.*, **37/38** (1992), 65–71.

[14] J. Cannon and W. Bosma, *CAYLEY. Quick Reference Guide*, Sydney, October 1991.

[15] G. E. Carlsson, J. E. Cruthirds, H. B. Sexton and C. G. Wright, "Interconnection networks based on a generalization of cube-connected cycles," *IEEE Trans. Computers*, **34** (1985), 769–777.

[16] D. Chudnovsky, G. Chudnovsky and M. Denneau, "Regular graphs with small diameter as models for interconnection networks," *3rd Int. Conf. on Supercomputing, Boston, International Computing Institute*, May 1988, 232–239.

[17] F. Comellas and J. Gómez, "New large graphs with given degree and diameter," *7th Internat. Conference on Graph Theory*, Kalamazoo, June 1992.

[18] R. Damerell, "On Moore graphs," *Proc. Cambridge Phil. Soc.*, **74** (1973), 227–236.

[19] M. J. Dinneen, *Algebraic Methods for Efficient Network Constructions*, Master's Thesis, Department of Computer Science, University of Victoria, Victoria, B. C., Canada, 1991.

[20] K. W. Doty, "New designs for dense processor interconnection networks," *IEEE Trans. Computers*, **C-33** (1984), 447–450.

[21] B. Elspas, "Topological constraints on interconnection-limited logic," *Proc. 5th Ann. Symp. Switching Circuit Theory and Logic Design* (1964),133–147.

[22] P. Erdős, S. Fajtlowicz and A. J. Hoffman, "Maximum degree in graphs of diameter 2," *Networks*, **10** (1980), 87–90.

[23] J. Gómez, M. A. Fiol and O. Serra, "On large $(\Delta, D)$–graphs," *Discrete Mathematics*, **114** (1993), 219–235.

[24] A. J. Hoffman and R. R. Singleton, "On Moore graphs with diameters 2 and 3," *IBM J. Res. Develop.*, **64** (1960), 15–21.

[25] G. Sabidussi, "Vertex transitive graphs," *Monatsh. Math.*, **68** (1969), 426–438.

[26] R. Storwick, "Improved construction techniques for $(d, k)$ graphs," *IEEE Trans. Computers*, **C-19** (1970), 1214–1216.